\documentclass[letterpaper]{article}
\usepackage[utf8]{inputenc}
\usepackage[T1]{fontenc}
\usepackage{hyperref}
\usepackage{url}
\usepackage{graphicx}
\usepackage{multicol}
\usepackage{subcaption}
\usepackage{tabularx}
\usepackage{amsfonts}
\usepackage{amsthm}
\usepackage{amssymb}
\usepackage{amsmath}
\usepackage{dsfont}
\usepackage{thmtools}
\usepackage{algorithm}
\usepackage{algpseudocode}
\usepackage[colorinlistoftodos,prependcaption,textsize=tiny]{todonotes}
\usepackage{lineno}
\usepackage{sidenotes}
\usepackage{marginnote}
\usepackage{array}
\usepackage{multirow}
\usepackage{authblk}
\usepackage{mathrsfs}
 
\usepackage[top=2cm, bottom=2cm, left=2cm, right=2cm]{geometry}

\title{Stable coexistence in indefinitely large systems of competing species}

\author[1,2,*]{M.N.~Mooij}
\author[2,3,4,5]{M.~Baudena}
\author[2,6]{A.S.~von der Heydt}
\author[1,2]{I.~Kryven}

\affil[1]{Mathematical Institute, Utrecht University, Utrecht, The Netherlands}
\affil[2]{Centre for Complex Systems Studies (CCSS), Utrecht University, Utrecht, The Netherlands}
\affil[3]{National Research Council of Italy, Institute of Atmospheric Sciences and Climate (CNR-ISAC), 10133 Torino, Italy}
\affil[4]{National Biodiversity Future Center (NBFC), 90133 Palermo, Italy}
\affil[5]{Copernicus Institute of Sustainable Development, Utrecht University, Utrecht, The Netherlands}
\affil[6]{Institute for Marine and Atmospheric Research (IMAU), Utrecht University, Utrecht, The Netherlands}

\affil[*]{m.n.mooij@uu.nl}

\date{}
\begin{document}

\maketitle

\section{Abstract}
The Lotka-Volterra system is a set of ordinary differential equations describing growth of interacting ecological species. 
This model has gained renewed interest in the context of random interaction networks.
One of the debated questions is understanding how the number of species in the system, $n$, influences the stability of the model.
Robert~May demonstrated that large systems become unstable, unless species-species interactions vanish.
This outcome has frequently been interpreted as a universal phenomenon and summarised as "large systems are unstable". 
However, May's results were performed on a specific type of graphs (Erd\H{o}s-R\'enyi), whereas we explore a different class of networks and we show that the competitive Lotka-Volterra system maintains stability even in the limit of large $n$, despite non-vanishing interaction strength.
We establish a lower bound on the interspecific interaction strength, formulated in terms of the maximum and minimum degrees of the ecological network, rather than being dependent upon the network's size.
For values below this threshold, coexistence of all species is attained in the asymptotic limit. In other words, the outlier nodes with large degree cause instability, rather than the large number of species in the system.
Our result refines May's bound, by showing that the type of network model is relevant and can lead to completely different results.
\\ \\
{\noindent \small \bf Keywords: Coexistence, Generalised Lotka-Volterra, Competition, Networks}

\section{Introduction}
Robert~May exposed the role of random network structures in the description of large ecological systems in his seminal work~\cite{may1972will}. 
A key conclusion drawn from his work is the diversity-stability paradox, which states that the likelihood of stable coexistence-defined as the convergence to a state where all species are present with positive density-is negatively affected by the size of the system $n$.
Thus, large networks tend to be less stable than smaller ones. 
This is a paradox since empirical observations indicate that large ecological networks are ubiquitous and stable~\cite{mccann2000diversity}. 
May studied the Erd\H{o}s R\'enyi model with edge weights and established that for random interaction matrices with sufficiently small average interaction, decaying at least as $n^{-1/2}$, coexistence is guaranteed with probability equal to one.
Consequently, large systems are less favourable to stable coexistence compared to small systems. 
However, since it is unclear whether this conclusion is universally applicable, a host of research followed to verify May's principle on networks with other connectivity assumptions, searching for structural patterns that stabilise large networks.

Allesina and Tang~\cite{allesina2012stability} and Allesina et al.~\cite{allesina2015} generalised May's model to include exclusively competitive, exclusively mutualistic, and predator-prey interaction types. 
Mutualistic interactions occur when two species mutually benefit each other, whereas competitive interactions arise when two species adversely affect each other.
Predator-prey interactions entail an asymmetric interaction, as observed, for instance, when a predator consumes its prey (henceforth the corresponding terminology).
The inverse relationship between system complexity (size) and stability has been demonstrated under all these interaction types. 

Due to the observed diversity-stability paradox in numerous systems, a substantial body of work has aimed to explain the presence of large systems in nature, using a variety of techniques and ideas.
Serv\'an et al.~\cite{servan2018} developed an evolutionary model that produces large stable networks by selectively removing "weak" species.
Bairey, Kelsic and Kishony~\cite{bairey2016} argued that while stability of networks governed solely by pairwise interactions is hampered by a large number of species, the inclusion of higher-order interactions may yield opposite effects, when mutualistic interactions are considered.
Grilli et al.~\cite{grilli2017} and Eppinga et al.~\cite{eppinga2018} considered higher-order interactions in competitive network models, leading to similar conclusions.
Grilli, Rogers, and Allesina~\cite{grilli2016} studied the impact of network compartmentalisation on stability.
Barabás, Michalska-Smith  and Allesina~\cite{barabas2017} showed that negative self-regulation has a favourable effect on stability.
Recently, Hatton et al.~\cite{hatton2024} have shown that using sublinear growth instead of logistic growth can account for diversity in ecosystems.
For an exhaustive review on complexity (number of species) and stability, we refer the reader to Landi et al.~\cite{landi2018}.

In this study, we explore a bounded degree network topology, wherein the maximum degree of the network is constrained by a function independent of the number of species.
This contrasts May's assumptions, which allowed nodes of arbitrary large degrees.
Our assumption of bounded degree, often referred to as network sparseness, is ecologically reasonable and also motivated by physical considerations. 
For instance, ecological networks inherently possess an underlying spatial domain with locally finite carrying capacity, possibly limiting the number of interactions in each location.
Moreover, ecological networks are likely to be sparse for structural reasons, as detailed for mutualistic networks by refs.~\cite{jordano2016sampling, bascompte2013mutualistic}. In fact, many interactions might be forbidden by different natural historic reasons (e.g. asynchronicities connected to phenology, especially in highly seasonal environments).
Sparseness in recorded ecosystem data may also be overestimated due to the inability to measure all interactions.
This limitation is particularly pronounced in the case of competitive interactions, as they are an indirect effect exerted by a third party, such as a shared food resource or limited habitat space. 
Nonetheless, undersampling of interactions is not believed to be the major reason behind unobserved links~\cite{jordano2016sampling}. In tropical forests, plant-plant heterospecific interactions have been observed to be very often close to zero~\cite{comita2010}, pointing at the fact that the interaction matrix might be sparse in these very relevant and extremely diverse systems. Moreover, sparse species interactions have been identified as a distinctive structural property in microbial communities~\cite{camacho2024}.

We will use the generalised Lotka-Volterra equations (well known in ecology but also celebrated in for example economics~\cite{comes2012banking, moran2019may}), which characterise the system in terms of the growth rates $\mathbf{r}$ and the symmetric interaction matrix $M$, are expressed as follows:
\begin{align}\label{eq:generalsystem}
    \frac{\mathrm{d} \mathbf{x(t)}}{\mathrm{d} t} = \mathbf{x(t)} \circ ( \mathbf{r} + M \mathbf{x(t)} ), 
\end{align}
where $\mathbf{x(t)}$ represents the species abundances at time $t \geq 0$, and $\circ$ denotes the Hadamard (component-wise) product. Without loss of generality, we impose the initial condition $\mathbf{x(0)}={\mathbf x}_0\in[0,\infty)^n$.

Different types of ecological interactions are modelled by choosing appropriate pair-wise interaction matrix $M$. 
Positive values for the elements $M_{ij}$ and $M_{ji}$ indicate mutualism between species $i$ and species $j$. Negative values of $M_{ij}$ and $M_{ji}$ indicate a competitive relationship, whereas positive $M_{ij}$ and negative $M_{ji}$ imply a predator-prey dynamic.
\begin{figure}[h!]
    \centering
    \includegraphics[width=0.8\linewidth]{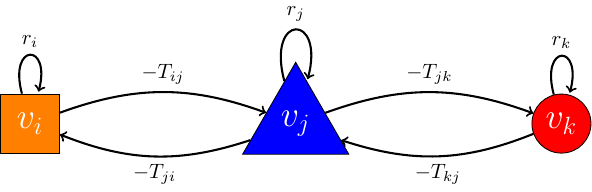}
    \caption{The conceptual structure of a three-species network ($v_i$, $v_j$, and $v_k$) within a competitive ecological system. The interactions $\it T$ (black arrows between species) are symmetric and have negative sign, while the growth rates $\mathbf{r}$ (black arrows within one species) are all strictly positive.}
    \label{fig:generalinteractionecologicalfigure}
\end{figure}
In the rest of the paper we consider two types of competitive interactions on a network with bounded degrees.
In Section~\ref{sec:constant}, we consider a competitive system wherein each species has unit growth rate and the same (negative) interaction with every other species in its neighbourhood.
This is achieved by setting the growth rates as $\mathbf{r} = \mathbf{1}$ and the interaction matrix as $M = - \tau A - \mathbb{I}$, where $\tau > 0$ is a parameter representing the competition strength and $A$ represents the binary adjacency matrix of an undirected network with bounded degree.
Then, in Section~\ref{sec:general}, we examine an arbitrary competitive system, with growth rates defined by a positive vector $\mathbf{r}$, and an interaction matrix $M$ represented as $M = -T - D$. The matrix $T$ contains only non-negative entries, while the diagonal self-regulation matrix $D$ exclusively contains positive elements.
Note that the choice $\mathbf{r} = \mathbf{1}$, $T = \tau A$, and $D = \mathbb{I}$ yields the constant interaction competitive system, studied in the first case.
In both scenarios, we show the existence of a subset within the parameter space (specified in terms of $\mathbf{r}$, $T$, and $D$) where stable coexistence of the dynamics is achieved in the asymptotic limit. This contradicts previous findings where the size of the subset of the parameter space facilitating stable coexistence was observed to shrink as $n^{-1/2}$.

\section{Results}
We define the coexistence of species as stable when, for certain interaction types (represented by $M$) and growth rates $\bf r$, and every initial condition with positive components, system~\eqref{eq:generalsystem} converges to a steady-state solution with only positive components.
Generally speaking, system~\eqref{eq:generalsystem} can have up to $2^n$ steady-state solutions. However, among these solutions, only those that solve the condition $M\mathbf{x^*} = -\mathbf{r}$ can be component-wise positive.
To demonstrate stable coexistence, we will verify two conditions: i) Feasibility, ensuring that $\mathbf{x^*}$ is indeed component-wise positive, and ii) Stability, which requires that the latter steady state solution is globally stable. 
Goh~\cite{GohBS} demonstrated that the former condition (feasibility) implies the latter condition (stability) if symmetric 
$M$ is negative definite.
For assessing local stability, it is sufficient to study the Jacobian evaluated at $\mathbf{x^*}$, which is given by:
\begin{align}
    J(\mathbf{x^*})_{ij} &= x^*_{i} M_{ij}. \label{eq:jacobian}
\end{align}

\subsection{Case 1: constant competition}\label{sec:constant}
\begin{figure}[h!]
    \centering
    \includegraphics[width=1\linewidth]{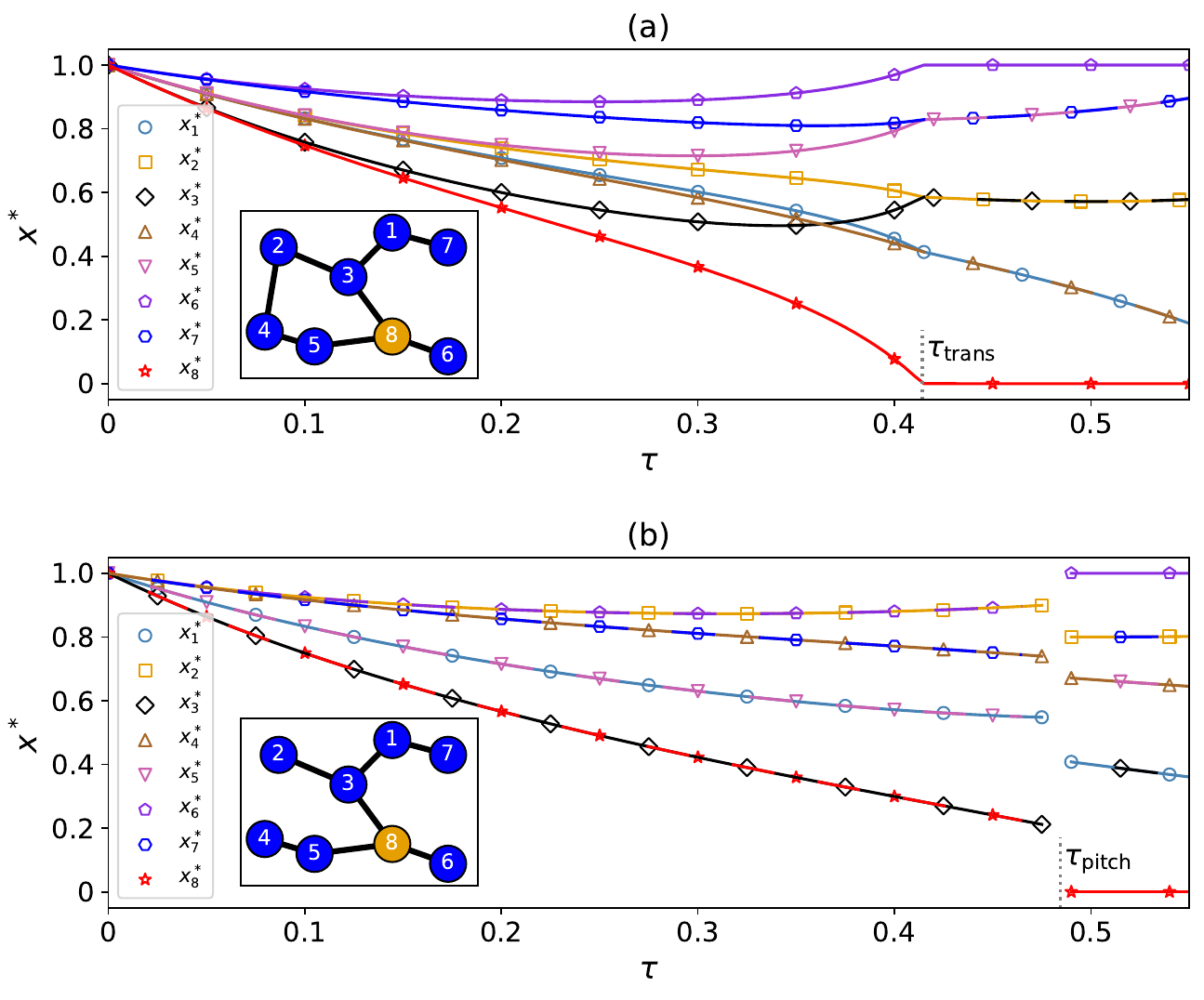}
    \caption{
    Two possible bifurcation types for constant interaction competitive Lotka-Volterra systems with different underlying interaction networks.
    In the networks, we highlighted the vertex corresponding to the variable going to zero, identified as vertex $8$.
    Note that the lower network is derived from the upper one by simply removing the edge $(2, 4)$.
    {\bf (a)} In this network a transcritical bifurcation occurs, where the bifurcation value of the interaction parameter is denoted by $\tau_{trans}$. 
    {\bf (b)} In this (only slightly modified) network, a pitchfork bifurcation takes place, demonstrating a non-continuous phase transition at the bifurcation value $\tau_{pitch}$. 
    Different variables are represented by different colors, with dashed lines featuring multiple colors to denote equal variable values.
    }
    \label{fig:bifurcationdiagram}
\end{figure}

In this scenario, the  interaction matrix $M$ is solely determined by the interaction strength parameter $\tau$, with stability ensured when $\tau=0$, and resulting in $\mathbf{x^*} = \mathbf{1}$.
Investigating $\mathbf{x^*(\tau)}$ reveals that there are two pathways, illustrated in Figure~\ref{fig:bifurcationdiagram}, where a discontinuity emerges upon incrementally increasing $\tau$, leading to a transition to a distinct steady state:
\begin{enumerate}
    \item a transcritical bifurcation, when $\mathbf{x^*}$ moves outside of the unit hypercube $[0,1]^n$ and therefore becomes unfeasible. This bifurcation is induced by the $\mathbf{x^*}$ term in Eq.~\eqref{eq:jacobian};
    \item a pitchfork bifurcation, when $\mathbf{x^*} \in (0,1)^n$ becomes locally unstable. This bifurcation is induced by the $M$ term in Eq.~\eqref{eq:jacobian};
\end{enumerate}

\noindent Note that the transcritical bifurcation entails a continuously declining coexistence point, whereas a pitchfork bifurcation occurs abruptly and non-continuously. 
For both types of bifurcation, we will show a lower bound for the parameter value at which the bifurcation occurs; first for the transcritical bifurcation, and secondly for the pitchfork bifurcation.

Rewriting the $k$th component of the stable coexistence point $\mathbf{x^*} = -M^{-1}\mathbf{r}$ as a Neumann series yields (possible for $\tau < d_{\text{max}}^{-1}$):
\begin{align}
    x^*_k &= \sum_{i=0}^{\infty} (-1)^i \tau^i \left( A^i \mathbf{1} \right)_k  \\
      &\geq \sum_{l=0}^{\infty} \tau^{2l} d_{\text{min}}^{2l}  - \sum_{l=1}^{\infty} \tau^{2l-1} d_{\text{max}}^{2l-1} \\
    &= \frac{ 1 - d_{\text{max}} \tau - d_{\text{max}}^2 \tau^2 + d_{\text{min}}^2 d_{\text{max}} \tau^3 }{ ( d_{\text{min}}^2 \tau^2 - 1 )( d_{\text{max}}^2 \tau^2 - 1 ) }.
\end{align}

\noindent
For the inequality, we utilized that the powers of the adjacency matrix enumerate the number of distinct walks, which can be bounded using the minimum and maximum degrees, $d_{\text{min}}$ and $d_{\text{max}}$, of the network.
Let $\alpha = d_{\text{min}} / d_{\text{max}}$ and consider the expression
\begin{align}
    P(\tau) &:= d_{\text{max}}^3 \alpha^2 \tau^3 - d_{\text{max}}^2 \tau^2 - d_{\text{max}} \tau + 1    
\end{align}
The transcritical bifurcation does not occur as long as $P(\tau) > 0$; an analysis of the third-order polynomial equation $P(\tau)=0$ reveals that this equation has three real roots, with the smallest one being:
\begin{align} \label{eq:constantcompetitivesystem}
    \begin{split}
        \Omega &\equiv \frac{ 1 - 2 \sqrt{1 + 3 \alpha^2} }{ 3 \alpha^2 } \sin\bigg[ \frac{1}{3} \arcsin\left( \theta \right) \bigg] \frac{1}{d_{\text{max}}} \geq \frac{\sqrt{5}-1}{2} \frac{1}{d_{\text{max}}}, \\
        \theta &= \frac{2 + 9 \alpha^2 - 27 \alpha^4}{2 (1 + 3 \alpha^2)^{3/2}}.
    \end{split}
\end{align}
This expression provides an upper bound on $\tau$ for which we can ensure that no transcritical bifurcation has occured. 
It remains to show that the pitchfork bifurcation will not occur for $\tau$ values smaller than the bound given by Eq.~\ref{eq:constantcompetitivesystem}.
By Gershgorin's Circle Theorem, we know that all eigenvalues of the matrix $M$ are contained in the union of the disks $\mathscr{D}_{-1}\left(\sum\limits_{j \neq i} |A_{ij}|\right)$, were $\mathscr{D}_{p}(r)$ is the disk centered around the (possibly complex) point $p$, with radius $r>0$. 
In this scenario, this simplifies to:
\begin{align}
    \mathscr{D}_i &= \mathscr{D}_{-1}\left(\tau d_i\right),
\end{align}
where $d_i$ denotes the degree of vertex $v_i$. 
Hence, if the condition $\tau d_i < 1$ holds true for all $i = 1,2,\dots, n$, we can ascertain that the coexistence point is stable whenever it is feasible. 
This condition is indeed fulfilled for $\tau < \Omega \leq d_{\text{max}}^{-1}$.

\begin{figure}[h!]
    \centering
    \includegraphics[width=\linewidth]{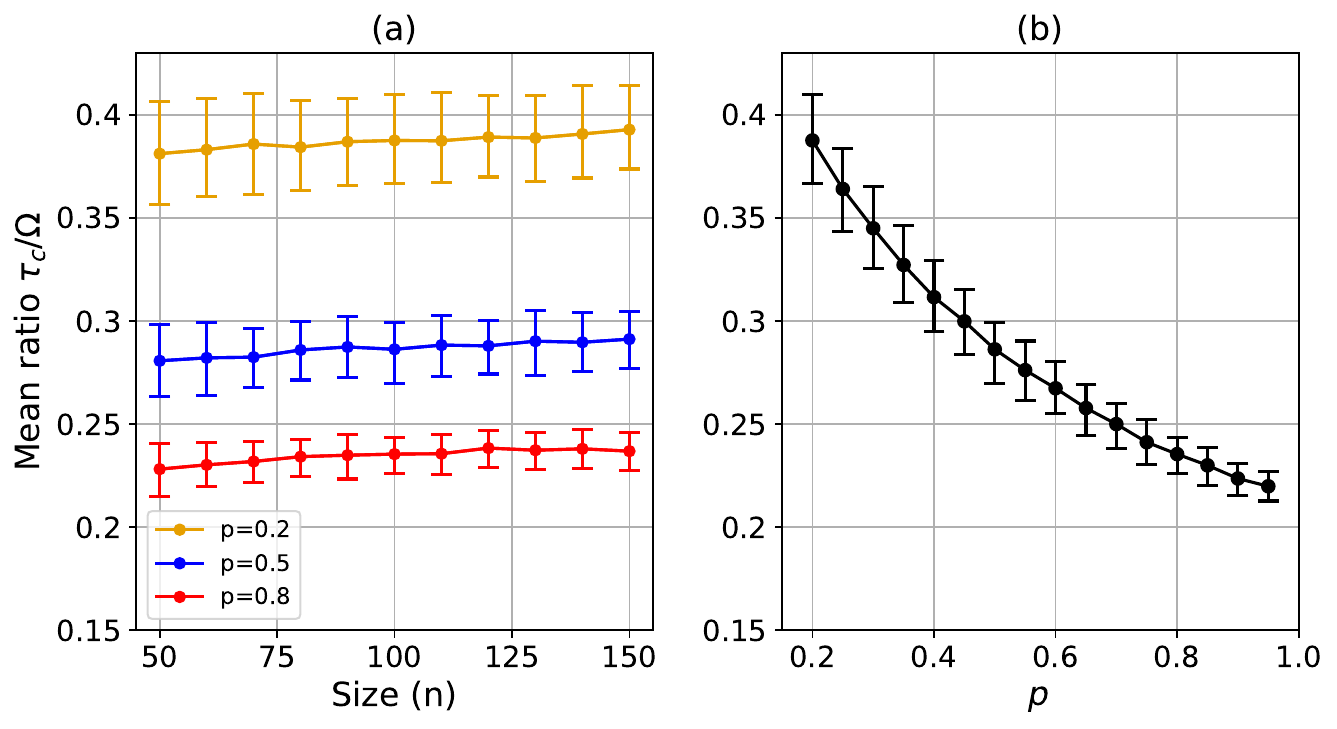}
    \caption{
    The ratio of the bifurcation value $\tau_{c}$ to our bound $\Omega$ as a function of the network size and connectivity. We generate networks with binomially distributed degree sequences using the parameters $d_{\text{max}}=30$, and varying the connectivity parameter $p$.
    Note that the maximum degree of any constructed network cannot exceed thirty.
    We verify whether the degree sequence is graphical, and if it is, we generate the network using the configuration model. 
    We use exhaustive search to ensure that the generated networks are connected.     
    {\bf(a)} The ratio of the bifurcation value $\tau_c$ to the conceived bound $\Omega$, averaged over $500$ runs. The error bars represent the 95\% confidence intervals of the data. Note the apparent independence of the ratio on the network size, while showing a strong dependency on the success probability $p$. 
    {\bf(b)} The ratio of the bifurcation value $\tau_c$ to the conceived bound $\Omega$ is determined for different connection probabilities $p$ and averaged over $500$ runs for each probability value. The error bars represent the 95\% confidence intervals of the data. The simulations were conducted on connected networks of size $100$, with binomially distributed degrees.    
    }
    \label{fig:binomialdegreessimulation}
\end{figure}

In conclusion, we have established that any interaction strength $\tau \in [0,\Omega)$ ensures stable species coexistence.
Note that the converse statement is not necessarily true; there may exist networks with $\tau > \Omega$ that still permit stable coexistence.
Remarkably, this interval is independent of the size of the network, indicating that coexistence is not contingent upon the network size when the maximum degree $d_{\text{max}}$ is fixed. 
If the maximum degree of the network is bounded, the probability of stable coexistence does not diminish to zero even as the system size grows indefinitely.
Our bound is tight across the entire spectrum of possible networks, as it aligns with the bifurcation point for complete bipartite graphs~\cite{MIS_paper}.
However, by restricting attention to a particular subclass of graphs, it is possible to determine the bifurcation value more precisely.
We demonstrate this with the following two examples.
\\ \\
{\bf Example 1. Regular graphs.} \\
In regular graphs all degrees are the same and the maximum eigenvalue of the adjacency matrix is given by the constant degree $d$, corresponding to the eigenvector $\mathbf{1}$. Therefore, we know
\begin{align}
    \mathbf{x^*} := -M^{-1}\mathbf{1} &= ( \tau A + \mathbb{I} )^{-1} \mathbf{1} = (\tau d + 1)^{-1} \mathbf{1}.
\end{align}
Let $\mathbf{v}\in \mathbb{R}^n$ be an eigenvector of $A$ corresponding to eigenvalue $\lambda$; then there exists an eigenvalue of the Jacobian evaluated at $\mathbf{x^*}$, given by
\begin{align}
    J(\mathbf{x^*}) \mathbf{v} &= \text{diag}(-M^{-1}\mathbf{1}) M \mathbf{v}  = -\frac{\tau \lambda + 1}{\tau d + 1} \mathbf{v}. \label{eq:regulargraphsjacobian}
\end{align}
Note that the stability of the interior fixed point relies on the topology of the regular network. 
In Eq.~\eqref{eq:regulargraphsjacobian}, we observe a monotonically increasing eigenvalue with respect to $\tau$. 
Hence, by using the Alon-Boppana bound, we can infer that almost surely (with probability one) every eigenvalue of the Jacobian adheres to the bound:
\begin{align}
    -1 \leq \lambda \leq - \frac{ -2 \tau \sqrt{d-1} + 1 }{ \tau d + 1 }.
\end{align}
Consequently, we obtain stable coexistence for almost all $d$-regular network configurations when $\tau < \frac{1}{ 2 \sqrt{d-1} }$. 
Note the scaling of order $d^{-1/2}$ compared to the $d^{-1}$ scaling for general graphs, suggesting that regular networks enhance coexistence.
\\ \\
{\bf Example 2. Configuration model.} \\
To illustrate our findings, we have generated multiple networks with binomially distributed degrees using the configuration model. 
We use a binomial degree distribution to control the maximum degree of the networks, which we set to be $30$. 
For various network sizes, we generated $100$ such networks and computed the average ratio of the actual bifurcation value $\tau_c$ to our bound $\Omega$~(Eq.~\eqref{eq:constantcompetitivesystem}), as depicted in Figure~\ref{fig:binomialdegreessimulation}. 
This ratio indicates the deviation of the true bifurcation point from our bound. 
We conducted different simulations for various connectivity values $0 < p < 1$, where large $p$ tends to yield networks with larger degrees on average, while a value of $p=0$ corresponds to a disconnected network. 
Figure~\hyperref[fig:binomialdegreessimulation]{\ref*{fig:binomialdegreessimulation}a} illustrates what appears to be an independent behavior of the calculated ratio with respect to the size of the network. 
Indeed, this observation aligns with our analytical result indicating that stable coexistence can occur for large networks whenever the maximum degree of the network is constrained.
As illustrated in Figure~\hyperref[fig:binomialdegreessimulation]{\ref*{fig:binomialdegreessimulation}b}, the tightness of the bound does depend on $p$ and improves when $p$ is smaller.

\subsection{Case 2: general competition}\label{sec:general}
The previous result can be expanded to a broader setting, where growth rates are represented by a positive vector $\mathbf{r} \in \mathbb{R}_{>0}$, and the interaction matrix is defined as $M = -T - D$. Here, $T$ is a symmetric matrix with strictly non-negative entries, and $D$ is a diagonal matrix with strictly positive entries. Let $D_\text{min}$ and $D_\text{max}$ represent the minimum and maximum elements of the vector $D$, respectively, while $r_\text{min}$ and $r_\text{max}$ denote the minimum and maximum elements of the vector $r$. Let $\Delta$ denote the maximum outgoing edge weight of a vertex in $G$, and let $\beta$ represent the ratio between the minimum outgoing edge weight and $\Delta$.
Applying analogous steps as in Case 1, we derive an upper bound for the outgoing edge strength $\Delta$ below which a transcritical bifurcation will not take place:
\begin{align}
    x^{\ast}_k &:= \left( (T+D)^{-1}\mathbf{r} \right)_k \\
    &= \sum_{i=0}^{\infty} (-1)^i \bigg[ (D^{-1} T)^i D^{-1} \mathbf{r} \bigg]_k \\
    &\geq \frac{ D_{\text{max}} D_{\text{min}}^2 r_{\text{min}} - D_{\text{max}} r_{\text{min}} \Delta^2 - D_{\text{max}}^2 r_{\text{max}} \Delta + r_{\text{max}} \beta^2 \Delta^3 }{ ( D_{\text{max}}^2 - \beta^2 \Delta^2 ) ( D_{\text{min}}^2 - \Delta^2 ) }. \label{eq:polynomial_general}
\end{align}
Setting the numerator of the last fraction to zero and solving for $\Delta$, we obtain the bound:
\begin{align}
    \begin{split}
    \Omega &\equiv \frac{ D_{\text{max}}^2 - 2 \sqrt{D_{\text{max}}^4 + \frac{3 \beta^2 D_{\text{max}} r_{\text{min}}}{r_{\text{max}}}} \sin\bigg[ \frac{1}{3} \arcsin\left( \theta \right) \bigg] }{ 3 \beta^2 }, \\
    \theta &= \frac{ \sqrt{\frac{r_{\text{max}}}{D_{\text{max}}}} ( 2 D_{\text{max}}^5 r_{\text{max}} + 9 \beta^2 r_{\text{min}}(D_{\text{max}}^2 - 3 \beta^2 D_{\text{min}}^2) ) }{ 2( D_{\text{max}}^3 r_{\text{max}} + 3 \beta^2 r_{\text{min}} )^{3/2} },
    \end{split} \label{eq:generalcompetitivesystem}
\end{align}

\noindent To show that a pitchfork bifurcation does not occur for $\Delta<\Omega$, we constrain the eigenvalue spectrum by examining the disks in the complex plane centred around the points $-D_{ii}$, each with a radius $\sum\limits_{j \neq i} T_{ij}$. 
All these disks are situated within the negative complex half-plane when $\Delta := \max\limits_{i \in [n] } \sum\limits_{j \neq i} |T_{ij}| < D_{\text{min}}$. 
This establishes a bound on the maximum outgoing edge weight, ensuring that for all values below this threshold, a pitchfork bifurcation does not occur. 
Note that the numerator from Eq.~\ref{eq:polynomial_general}, evaluated at the point $D_\text{min}$, evaluates to $-D_{\text{max}}^2 + \beta^2 D_{\text{min}}^2$, which is always negative. 
The numerator of Eq.~\ref{eq:polynomial_general} at $D_\text{min}$ is $-D_{\text{max}}^2 + \beta^2 D_{\text{min}}^2$, which is always negative.
Hence, the bound $\Omega$ is always smaller than the bound for a pitchfork bifurcation.
Combining these observations, we establish an upper bound $\Omega$ on the maximum outgoing edge weight $\Delta$, thereby ensuring the absence of both transcritical and pitchfork bifurcations.

Similar to constant interaction systems, we observe the existence of a subset within the parameter space where stable coexistence is maintained.
Note that the derived bound does not diminish when the network size increases indefinitely, given fixed values for $D_\text{max}$ and $r_\text{max}$.
Hence, even in the limit of large $n$, stable coexistence persists under these parameters.

One can use the parameter set $\mathbf{r} = \mathbf{1}$, $D = \mathbb{I}$ and $T = \tau A$ to get back the original bound for the constant interaction case (Eq.~\eqref{eq:constantcompetitivesystem}). We examine different regimes of the parameter space, identified with ecologically meaningful cases, which are outlined in Table~\ref{tab:bound_regimes}. Large self-inhibition $D_{\text{max}}$ makes the system collapse, as the species die out too soon to have meaningful dynamics. Large maximum growth rate $r_{\text{max}}$ does not allow for coexistence, as the neighbours of $v_i$ are heavily influenced by node $v_i$, and therefore it does not allow for coexistence. In the case of vanishing $r_{\text{min}}$, we will also not obtain coexistence, as the corresponding species will not grow and therefore will die out in the long run.

\begin{table}[h]
\centering
\begin{tabular}{|c|c|c|c|}
\hline
regime & parameters & $\theta$ & $\Omega$ \\
\hline
constant interaction & $\mathbf{r}=\mathbf{1}$, $T=\tau A$, $D=\mathbb{I}$ & $\frac{2 + 9 \alpha^2 - 27 \alpha^4}{2 (1 + 3 \alpha^2)^{3/2}}$ & $\frac{ 1 - 2 \sqrt{1 + 3 \alpha^2} }{ 3 \alpha^2 } \sin\bigg[ \frac{1}{3} \arcsin\left( \theta \right) \bigg] \frac{1}{d_{\text{max}}}$ \\

large self-inhibition & $D_{\text{max}} \to \infty$ & 1 & 0 \\
large growth rate & $r_{\text{max}} \to \infty$ & 1 & 0 \\
small growth rate & $r_{\text{min}} \to 0$ & 1 & 0 \\

\hline
\end{tabular}
\caption{Asymptotic regimes of $\Omega$.}
\label{tab:bound_regimes}
\end{table}

Finally, the converse statement to our main result is also true; when the maximum degree $k$ diverges with network size, at least one species will necessarily become extinct. To see why this is the case, consider the star graph $S_{k}$, in which all but one node have degree $1$, and one node has degree $k-1$. 
Note that any graph $G$ with maximum degree $k$, contains $S_{k}$ as a subgraph.
The minimum eigenvalue of $S_{k}$ is $-(k-1)^{1/2}$, and due to the eigenvalue interlacing theorem, we know that $\lambda_{\text{min}}(A) \leq \lambda_{\text{min}}(S_{k})$. 
Consequently, (at least) one pitchfork bifurcation must occur for $\tau > (k-1)^{-1/2}$. 
Since $k$ grows indefinitely, this value converges to zero, indicating that a pitchfork bifurcation must occur for networks with diverging degree. 
This proof easily generalises to our general case of weighted graphs, as we can bound the minimum eigenvalue of the interaction matrix (due to negative definiteness) by the minimum eigenvalue of the binary adjacency matrix of the network scaled by the maximum element of $T$. 
Therefore, by the same reasoning, a pitchfork bifurcation must occur when $k$ diverges.

\section{Conclusions}
This study examines coexistence of competitive species in the Lotka-Volterra system with undirected network of symmetric competitive interactions. We prove that when the degrees are bounded, regardless of system size, there is a globally-stable coexistence state for a parameter set of nonzero measure.
This observation resolves the diversity-stability paradox by demonstrating that stable species coexistence is possible even in indefinitely large systems. When applied to weighted networks, we show that stable species coexistence can be maintained indefinitely, when the sum of edge weights is bounded for each vertex. The converse statement is also shown to be true: if the largest degree diverges with system size, this will necessarily lead to extinction of at least one species.

Additionally to resolving the paradox, an immediate implication of our results for the ecological modelling community is that a subtle choice between the network models that does or does not allow for limiting the maximum degree will have drastic implications on the dynamics of the competitively interacting species.
Rephrasing this in the language of statistical physics, choosing the interaction network from microcanonical or canonical ensembles, which are typically thought to be equivalent in other contexts, might determine the fate of species coexistence in the Lotka-Volterra model. 

It still remains an open question whether similar type of behaviour is also characteristic for directed networks and networks with mixed mutualistic-competitive interactions.
While interactions of this type have been explored in Erd\H{o}s-R\'enyi networks, their implications in bounded degree networks can lead to significantly different conclusions.
Even though we show that the diversity-stability paradox does not exist per se, we note that our finding does not diminish the relevance of other possible explanations of diversity, for example related to competition for time~\cite{levine2024} or space~\cite{tilman1994} or to higher order interactions~\cite{grilli2017, eppinga2018}.

Finally, let us mention that this result is potentially relevant in ecology, since networks with sparse structures, and hence presumably bounded degrees, are likely to be common in ecology~\cite{jordano2016sampling}, for example in the context of tropical forests~\cite{comita2010}, where the underlying interaction network is possibly better represented as sparse rather than as an Erd\H{o}s-R\'enyi graph.
When there is a rationale to assume a bounded degree network, it is essential to acknowledge that the conclusions drawn for other network models may no longer be applicable.

\section{Data availability}
All data that support the figures within this paper and other findings of
this study are available at\\\url{https://figshare.com/s/9dff23b8f6a6fc1a77ee}.

\section{Code availability}
Code is available for this paper at\\
\begin{minipage}{\linewidth}
  \href{https://github.com/NiekMooij/Stable-coexistence-in-indefinitely-large-systems-of-competing-species.git}{%
    \parbox{\linewidth}{%
      https://github.com/NiekMooij/Stable-coexistence-in-indefinitely-large-systems-of-competing-species.git}%
    }
\end{minipage}

\section{Acknowledgements}
NM gratefully acknowledges support from Complex Systems Fund, with special thanks to Peter Koeze. MB acknowledges the Italian National Biodiversity Future Center (NBFC): National Recovery and Resilience Plan (NRRP), Mission 4 Component 2 Investment 1.4 of the Italian Ministry of University and Research; funded by the EU - NextGenerationEU (Project code CN 00000033). The work of A.S.vdH was also funded by the Dutch Research Council (NWO) through the NWO-Vici project  ‘Interacting climate tipping elements: When does tipping cause tipping?’ (project VI.C.202.081). IK gratefully acknowledges support form Netherlands Research Organisation (NWO), research program VIDI, project number VI.Vidi.213.108. 
The authors would like to thank Maarten Eppinga for discussion and providing very helpful literature.

\section{Author contributions}
NM, MB, AvdH, IK designed the study. NM, IK performed the analytical analysis. NM developed the code, ran simulations, and made the figures. NM wrote the first draft, after which all authors contributed to improving the manuscript. MB, AvdH, IK provided supervision, with IK being the main supervisor.

\section{Competing interests}
None declared.

\section{Materials \& Correspondence}
Please contact NM for correspondence.

\end{document}